\documentclass[12pt,oneside,reqno]{amsart}
\usepackage{graphicx}
\usepackage{mathrsfs}
\usepackage{stmaryrd}
\usepackage{amsfonts}
\usepackage{amsmath}
\usepackage{enumerate,amsmath,amssymb,amsthm}
\newcommand{\be}{\begin{eqnarray}}
\newcommand{\ee}{\end{eqnarray}}
\newcommand{\ce}{\begin{eqnarray*}}
\newcommand{\de}{\end{eqnarray*}}
\newtheorem{theorem}{Theorem}[section]
\newtheorem{lemma}[theorem]{Lemma}
\newtheorem{remark}[theorem]{Remark}
\newtheorem{definition}[theorem]{Definition}
\newtheorem{proposition}[theorem]{Proposition}

\newtheorem{corollary}[theorem]{Corollary}
\def\HS{{\mathrm{\tiny HS}}}
\def\e{\varepsilon}

\def\a{\alpha}

\def\Om{\Omega}
\def\b{\beta}

\def\l{\lambda}

\def\[{{\Big[}}
\def\]{{\Big]}}
\def\<{{\langle}}
\def\>{{\rangle}}
\def\({{\Big(}}
\def\){{\Big)}}

\def\dif{{\mathord{{\rm d}}}}

\def\no{\nonumber}
\def\bt{\begin{theorem}}
\def\et{\end{theorem}}
\def\bl{\begin{lemma}}
\def\el{\end{lemma}}
\def\br{\begin{remark}}
\def\er{\end{remark}}
\def\bd{\begin{definition}}
\def\ed{\end{definition}}
\def\bp{\begin{proposition}}
\def\ep{\end{proposition}}
\def\bc{\begin{corollary}}
\def\ec{\end{corollary}}

\def\cF{{\mathcal F}}

\def\mE{{\mathbb E}}

\def\mR{{\mathbb R}}

\def\bP{{\mathbf P}}

\def\geq{\geqslant}
\def\leq{\leqslant}

\def\bE{{\mathbf E}}
\allowdisplaybreaks

\begin{document}
\title{Exponential Ergodicity of Non-Lipschitz Multivalued Stochastic Differential Equations$^{*}$}

\date{}
\author{Jiagang Ren$^{1}$, Jing Wu$^{1}$, Xicheng Zhang$^{2,3}$ }


\date{}

\thanks{$^{*}$ Research supported by NSFC (Grant no. 10871215)}



\maketitle

\noindent{$^{1}$ Zhongshan University, Guangzhou, Guangdong 510275, P.R.China\\
$^{2}$ The University of New South Wales, Sydney, 2052, Australia,\\
$^{3}$ Huazhong University of Science and Technology, Wuhan, Hubei 430074, P.R.China\\
Emails: J. Ren: renjg@mail.sysu.edu.cn; J. Wu:  wjjosie@hotmail.com;\\ X. Zhang: XichengZhang@gmail.com}
\vskip 0.4cm
{\bf Abstract}
\vskip 0.2cm
\begin{minipage}{15 cm}
We prove the exponential ergodicity of the transition probabilities of
solutions to elliptic multivalued stochastic differential equations.
\vskip 0.2cm
{\bf R\'esum\'e}
\vskip 0.2cm
On prouve l'ergodicit\'e exponentielle des
probabilit\'es de transition des equations diff\'erentielles stochastiques elliptiques.
\end{minipage}
\vskip 0.5cm
\section{Introduction and Preliminaries}
Consider the following stochastic differential equation: \be \label{Eq0} \dif X_t=b(X_t)\dif
t+\sigma(X_t)\dif W_t,\ \ X_0=x, \ee
where $b:\mR^d\to \mR^d$ and
$\sigma: \mR^d\to \mR^d\otimes \mR^n$ are continuous functions, $(W_t)_{t\geq 0}$ is an
$n$-dimensional standard Brownian motion defined on some complete
probability space $(\Om,\cF,\bP)$. When $\sigma$ is a uniformly
elliptic square matrix and $\sigma$ and $b$ satisfy some regular conditions (more precisely,
(H1), (H2) and (H4) below), it is recently proved in \cite{Zh} that the solution
is exponentially ergodic.

On the other hand, under the same uniform elliptic assumption and an
additional one that $\sigma$ and $b$ are $C_b^2$, C\'epa and Jacquot
proved in \cite{ce3} the ergodicity for the solution of the
following stochastic variational inequality (SVI in short): \be
\label{Eq1} \dif X_t+\partial\varphi(X_t)\ni b(X_t)\dif
t+\sigma(X_t)\dif W_t,\ \ X_0=x\in\overline{\mathrm{Dom}(\varphi)},
\ee where $\partial\varphi$ is the sub-differential of some convex
function $\varphi$ with a compact domain
$\mathrm{Dom}(\varphi)=\{x:\varphi(x)<\infty\}$.

A common drawback of the above two papers is the {\it uniform
elliptic assumption} of the diffusion coefficients. The purpose of
the present paper is to remove this assumption and instead assume
only the {\it ellipticity}. Our main result as stated in Theorem
\ref{Th3} below unifies and improves the main results of both of
\cite{Zh} and \cite{ce3}. In particular, our result applies to
stochastic variational inequalities defined on non-compact domains.
 Furthermore, we do not need to assume that the
diffusion matrix is square  and  our method even works for general
multivalued stochastic differential equations (MSDEs in
abbreviation): \be \label{Eq2} \dif X_t+A(X_t)\ni b(X_t)\dif
t+\sigma(X_t)\dif W_t,\ \ X_0=x\in\overline{D(A)}, \ee where $A$
is a multivalued maximal monotone operator on $\mR^d$ with
Int$(D(A))\neq\emptyset$.

Now we introduce notions and notations.
Given an operator $A$ from $\mR^d$
to $2^{\mR^d}$, define: \ce
D(A)&:=&\{x\in \mR^d:A(x)\neq \emptyset\}, \\
\mathrm{Gr}(A)&:=& \{(x,y)\in \mR^{2d}: x\in \mR^d, y\in A(x)\}. \de
Then $A$ is called monotone if $
\<y_1-y_2,x_1-x_2\>\geq 0$ for any  $(x_1,y_1),(x_2,y_2)\in
\mathrm{Gr}(A)$, and  $A$ is called maximal monotone if
\ce (x_1,y_1)\in \mathrm{Gr}(A)\Leftrightarrow\<y_1-y_2,x_1-x_2\>
\geq 0,\quad\forall(x_2,y_2)\in \mathrm{Gr}(A).\de

 \bd\label{d2}
A pair of continuous and $(\cF_t)$-adapted processes $(X,K)$ is called
a solution of (\ref{Eq2}) if

(i)  $X_0=x_0$, $X_t\in \overline{D(A)}$ a.s.;

(ii) $K$ is of locally finite variation and $K_0=0$ a.s.;

(iii) $\dif X_t=b(t,X_t)\dif t+\sigma(t,X_t)\dif W_t-\dif K_t$, $0\leq t<\infty,
\quad a.s.$;

(iv) For any continuous and $(\cF_t)-$adapted functions $(\a,\b)$ with
$(\a_t,\b_t)\in \mathrm{Gr}(A),\quad \forall t\in[0,+\infty)$,
the measure $\<X_t-\a_t, \dif K_t-\b_t\dif t\>$ is positive.
\ed

We make the following assumptions:
\begin{enumerate}[\bf ($\mathbf{H}1$)]

\item (Monotonicity) There exists $\lambda_0\in\mR$ such that for all
$x,y\in\mR^d$
$$
2\<x-y,b(x)-b(y)\>+\|\sigma(x)-\sigma(y)\|^2_{\HS}\leq\l_0|x-y|^2(1\vee\log |x-y|^{-1}).
$$

\item (Growth of $\sigma$) There exists  $\lambda_1>0$ such that for all $x\in\mR^d$
$$
\|\sigma(x)\|_\HS\leq \l_1(1+|x|).
$$

\item (Ellipticity of $\sigma$)
$$
\sigma\sigma^*(x)>0,~~~~~~~~~~~~~~\forall x\in \mR^d.
$$
\item (One side growth of $b$) There exist a $p\geq 2$ and constants  $\l_3>0,\l_4\geq 0$
such that for all $x\in\mR^d$
$$
2\<x,b(x)\>+\|\sigma(x)\|^2_\HS\leq -\l_3|x|^{p}+\l_4.
$$
\end{enumerate}

\bt\label{Th2} Assume {\bf($\mathbf{H}1$)} and {\bf($\mathbf{H}2$)}
hold. Then (\ref{Eq2}) has a unique strong solution.\et

\begin{proof}
 The existence of a weak solution is proved in
\cite {ce3} and the pathwise uniqueness can be proved in a more or less
standard way using a version of Bihari inequality (see \cite{rz1}).
Finally by Yamada-Watanabe's theorem the existence of a unique
strong solution follows.
\end{proof}

Let $\{X_t(x), t\geq 0,
x\in\mE\}$ denote the unique solution to
(\ref{Eq2}).
It is obviously a Markow family and its transition semigroup
and transition probability are defined
respectively as:
$$
P_tf(x_0):=\bE f(X_t(x_0)),~~t>0,~~f\in B_b(\mR^d)$$
and
$$
P_t(x_0,E):=\bP(X_t(x_0)\in E),
$$
where $x_0\in\mE$ and $B_b(\mR^d)$ denotes the set of all bounded
measurable functions on $\mR^d$. For general notions (e.g., strong
Feller property, irreducibility, ergodicity, etc)
 concerning Markov semigroups, we refer to \cite{ce3, Zh}.

\section{Main Result}
Now we state the main result of the paper.
\bt\label{Th3} Assume {\bf($\mathbf{H}1$)}-{\bf($\mathbf{H}3$)}.
Then the transition probability $P_t$ of the solution to (\ref{Eq2})
is irreducible and strong Feller. If in addition,
{\bf($\mathbf{H}4$)} holds, then there exists a unique invariant
probability measure $\mu$ of $P_t$ having full support in
$\overline{D(A)}$ such that
\begin{enumerate}[(i)]
\item If $p\geq2$ in {\bf($\mathbf{H}4$)}, then
for all $t>0$ and $x_0\in\overline{D(A)}$, $\mu$ is equivalent to
$P_t(x_0,\cdot)$, and
$$
\lim_{t\rightarrow\infty}\|P_t(x_0,\cdot)-\mu\|_{\mathrm{Var}}=0,
$$
where $\|\cdot\|_{\mathrm{Var}}$ denotes the total variation of a
signed measure.
\item If $p>2$ in {\bf($\mathbf{H}4$)}, then for  some $\a,C>0$
independent of $x_0$ and $t$,
$$
\|P_t(x_0,\cdot)-\mu\|_{\mathrm{Var}}\leq C\cdot e^{-\a t}.
$$
Moreover, for any $q>1$ and each $\varphi\in
L^q(\overline{D(A)},\mu)$
$$
\|P_t\varphi-\mu(\varphi)\|_{q}\leq C_q\cdot e^{-\a
t/q}\|\varphi\|_{q},\  \ \forall t>0,
$$
where $\a$ is the same as above and
$\mu(\varphi):=\int_{\overline{D(A)}}\varphi(x)\mu(\dif x)$. In
particular, let $L_q$ be the generator of $P_t$ in
$L^q(\overline{D(A)},\mu)$. Then $L_q$ has a spectral gap (greater
than $\a/q$) in $L^q(\overline{D(A)},\mu)$.
\end{enumerate}\et
The proof consists in proving the irreducibility and strong Feller property.
\subsection{Irreducibility}
\bl\label{l3} Suppose $y_0\in$ Int$(D(A)),$ $m>0$, and  $Y_t$ is the
solution to the following MSDE: \ce\label{Eq5}\dif Y_t+A(Y_t)\dif t\ni
-m(Y_t-y_0)\dif t+\sigma(Y_t)\dif W_t,~~Y_0=x_0, \de where $\sigma$ is
the diffusion coefficient of (\ref{Eq2}). Then under
{\bf($\mathbf{H}1$)} and {\bf($\mathbf{H}2$)} we have \ce
\bE|Y_t-y_0|^2\leq e^{-C(m)t}|x_0-y_0|^2+\frac{C_0}{C(m)},\de where
$C(m)=2(m-2\l_1^2-1/2)$ and
$C_0=2\l_1^2(1+2|y_0|^2)+|A^{\circ}(y_0)|^2.$
Here $A^{\circ}$ is the minimal section of $A$ and $|A^\circ(y_0)|<+\infty$ because $y_0\in$ Int$(D(A))$
(see \cite{ce1}).\el

\begin{proof}
The proof is adapted from \cite{ce3}. Consider the solution $Y_t^n$ to the following equation:
\ce \dif
Y_t^n+A_n(Y_t^n)\dif t=-m(Y_t^n-y_0)\dif t+\sigma(Y_t^n)\dif
W_t,~~Y_0^n=x_0
\de
where $A_n$ is the Yosida approximation of $A$. From \cite{ce1}
we know that $A_n$ is monotone, single-valued and $|A_n(x)|\nearrow|A^\circ(x)|$
if $x\in D(A)$, where $A^\circ$ is the minimal section of $A$.
Moreover, since the law of $Y_t^n$ converges to that of
$Y_t$, it is  enough to prove the inequality for
$Y_t^n$. Hence by {\bf($\mathbf{H}2$)}
\ce
&&-2m|x-y_0|^2+\|\sigma(x)\|_\HS^2-2\<A_n(x),x-y_0\>\\
&\leq&-2m|x-y_0|^2+\l_1^2(1+|x|)^2-2\<A_n(x)-A_n(y_0),x-y_0\>-2\<A_n(y_0),x-y_0\>\\
&\leq&-2m|x-y_0|^2+\l_1^2(1+|x|)^2+|x-y_0|^2+|A^{\circ}(y_0)|^2\\
&\leq&-2m|x-y_0|^2+2\l_1^2(1+2|x-y_0|^2+2|y_0|^2)+|x-y_0|^2+|A^{\circ}(y_0)|^2\\
&=&-C(m)|x-y_0|^2+C_0, \de
Thus, by It\^o's formula we have
\ce
\frac{d}{\dif t}\bE|Y_t^n-y_0|^2
&=&-2\bE(\<Y_t^n-y_0,A_n(Y_t^n)\>)+\bE[\mathrm{Tr}(\sigma\sigma^{*}(Y_t^n))]-2m\bE|Y_t^n-y_0|^2\\
&\leq&-C(m)\bE|Y_t^n-y_0|^2 +C_0. \de
Therefore \ce \bE|Y_t^n-y_0|^2\leq
e^{-C(m)t}|x_0-y_0|^2+\frac{C_0}{C(m)}.\de
\end{proof}
 \bp\label{p4} Under
{\bf($\mathbf{H}1$)}-{\bf($\mathbf{H}3$)}, the transition
probability $P_t$ is irreducible.\ep
\begin{proof}
It suffices to prove that for any
$x_0\in \overline{D(A)}$, $T>0$, $y_0\in$ Int$(D(A))$ and $a>0,$
$$
P_T(x_0,B(y_0,a))=\bP(X_T(x_0)\in B(y_0,a))=\bP(|X_T(x_0)-y_0|\leq
a)>0,$$ or equivalently:
$$
\bP(|X_T(x_0)-y_0|> a)<1.$$

Fix $a$, $T$ and $y_0$. By Lemma \ref{l3} and Chebyshev's inequality,
we can choose an $m$ large enough such that, denoting by $(Y_t,\tilde{K}_t)$ the unique solution to
\be\label{Eq6}\dif Y_t+A(Y_t)\dif t\ni -m(Y_t-y_0)\dif
t+\sigma(Y_t)\dif W_t,~~Y_0=x_0\in \overline{D(A)}, \ee
we have
\be
\bP(|Y_T(x_0)-y_0|>a)\leq
\left(e^{-C(m)T}|x_0-y_0|^2+\frac{C_0}{C(m)}\right)/a^2<1. \ee

Set
$$
\tau_N:=\inf\{t: |Y_t|\geq N\}.
$$
Note that by \cite{ce3}
$$
\bE\left[\sup_{t\in[0,T]}|Y_t(x_0)|\right]\leq C
$$
for some constant $C$ depending on $x_0, y_0, \lambda_1, m$ and $T$.
Thus we may fix an $N$ so that
\be\label{qandp}
\bP(\tau_N\leq T)+\bP(|Y_T(x_0)-y_0|>a)<1. \ee
Define
$$
U_t:=\sigma(Y_t)^{*}[\sigma(Y_t)\sigma(Y_t)^{*}]^{-1}(-m(Y_t-y_0)-b(Y_t))
$$
and
$$Z_T=\exp\left(\int_0^{T\wedge\tau_N}U_s\dif W_s-
\frac{1}{2}\int_0^{T\wedge\tau_N}|U_s|^2\dif s\right).
$$
Since $|U_{t\wedge \tau_N}|^2$ is bounded, $\bE[Z_T]=1$ by Novikov's criteria.

By Girsanov's theorem, $W_t^*:=W_t+V_t$ is a $Q$-Brownian motion, where
$$
V_t:=\int_0^{t\wedge\tau_N}U_s\dif s,\quad Q:=Z_T\bP.
$$
By (\ref{qandp}) we have \be\label{qless1} Q(\{\tau_N\leq
T\}\cup\{|Y_T(x_0)-y_0|>a\})<1. \ee Note that the solution
$(Y_t,\tilde{K}_t)$ of (\ref{Eq6}) also solves the MSDE below
\ce
Y_t+\int_0^tA(Y_s)\dif s\ni\int_0^t\sigma(Y_s)\dif W_s^*
+\int_0^{t\wedge \tau_N}b(Y_s)\dif s-\int_{t\wedge \tau_N}^t
m(Y_s-y_0)\dif s. \de Set
$$
\theta_N:=\inf\{t: |X_t|\geq N\}.
$$
Then the uniqueness in distribution for (\ref{Eq6}) yields that the
law of $\{(X_t\textbf{1}_{\{\theta_N\geq T\}})_{t\in [0,T]},
\theta_N\}$ under $\bP$ is the same as that of
$\{(Y_t\textbf{1}_{\{\tau_N\geq T\}})_{t\in [0,T]},\tau_N\}$
under $Q$. Hence
\ce
\bP(|X_T(x_0)-y_0|>a)&\leq &\bP(\{\theta_N\leq T\}\cup\{\theta_N\geq T, |X_T(x_0)-y_0|>a\})\\
&=&Q(\{\tau_N\leq T\}\cup\{\tau_N\geq T, |Y_T(x_0)-y_0|>a\})\\
&\leq &Q(\{\tau_N\leq T\}\cup\{|Y_T(x_0)-y_0|>a\})<1. \de
\end{proof}

\subsection{Strong Feller Property} The proof of the following lemma is plain by using
Kolmogorov's lemma on path regularity of stochastic processes.
\bl\label{lemma31} Denote by $(X_t(x), K_t(x))$ the solution of
(\ref{Eq2}) with inital value $x$. Then for any $p>d,$ there exists
$t_p>0$ such that for all $r>0$
$$
\bE\left[\sup_{x\in D_r, s\leq t_p}|X_s(x)|^p\right]<\infty,
$$
where $D_r:=\overline{D(A)}\cap \{|x|\leq r\}$. \el

 \bp\label{p2}  Under {\bf($\mathbf{H}1$)}-{\bf($\mathbf{H}3$)},
 the semigroup $P_t$ is strong Feller. \ep

\begin{proof}We divide the proof into two steps.

Step 1: Assume that there exists a $\l_2>0$ such that
$\|[\sigma^*\sigma]^{-1}\|_\HS\leq \lambda_2$. Consider the following drift transformed MSDE:
\be
\label{Eq3}\left\{
\begin{array}{ll}
&\dif Y_t+A(Y_t)\dif t\ni b(Y_t)\dif t+\sigma(Y_t)\dif W_t+|x_0-y_0|^{\a}\frac{X_t-Y_t}{|X_t-Y_t|}\cdot 1_{\{X_t\not=Y_t\}}\cdot
1_{\{t<\tau\}}\dif t,\\
& Y_0~=y_0\in \overline{D(A)},
\end{array}
\right. \ee
where $\a\in(0,1)$, $X_t$ is the solution to (\ref{Eq2})
and $\tau$ is the coupling time given by
$$
\tau:=\inf\{t>0: |X_t-Y_t|=0\}.
$$
An argument similar to \cite{Zh} allows to prove it admits a unique solution.

For $T>0$ define
$$
U_T:=\exp\left[\int^{T\wedge\tau}_0\<\dif W_s,H(X_s,Y_s)\>
-\frac{1}{2}\int^{T\wedge\tau}_0|H(X_s,Y_s)|^2\dif s\right]
$$
and
$$
\tilde W_t:=W_t+\int^{t\wedge\tau}_0H(X_s,Y_s)\dif s,
$$
where
$$
H(x,y):=|x_0-y_0|^{\a}
\cdot\sigma^*(y)[\sigma\sigma^*(y)]^{-1}\frac{x-y}{|x-y|}.
$$
Since $\|[\sigma\sigma^*(y)]^{-1}\|_\HS\leq\lambda_2$, we have
$$
|H(x,y)|^2\leq \l_2\cdot|x_0-y_0|^{2\a}.
$$
Thus,
$$
\bE U_T=1\ \ \mbox{ and } \ \ \bE
U^2_T\leq\exp\left[\l_2T\cdot|x_0-y_0|^{2\a}\right].
$$
By the elementary inequality $e^r-1\leq r e^r$ for $r\geq 0 $,
we have for any $|x_0-y_0|\leq\eta$,
\be\label{Eq7}
(\bE|1-U_T|)^2&\leq& \bE|1-U_T|^2=\bE U_T^2-1\no\\
&\leq&\exp\left[\l_2T\cdot|x_0-y_0|^{2\a}\right]-1\no\\
&\leq& C_{T,\l_2,\eta}\cdot|x_0-y_0|^{2\a} \ee
and
\be\label{Eq8}
\left(\bE\left[(1+U_T)1_{\{\tau\geq T\}}\right]\right)^2&\leq& (3+\bE U_T^2)\cdot\bP(\tau\geq T)\no\\
&\leq&C_{T,\l_2,\eta} \cdot\bP((2T)\wedge\tau\geq T)\no\\
&\leq& C_{T,\l_2,\eta}\cdot\bE((2T)\wedge\tau)/T.
\ee
First applying It\^o's formula to $\sqrt{|Z_{t\wedge\tau}|^2+\e}$ where $Z_s:=X_s-Y_s$,
then letting $\e\downarrow 0$,  and finally taking
expectation, we have by {\bf($\mathbf{H}1$)},
\ce
\bE|X_{t\wedge\tau}-Y_{t\wedge\tau}|\leq |x_0-y_0|-|x_0-y_0|^{\a}\cdot \bE(t\wedge\tau)
+\frac{\l_0}{2}\int^t_0\rho_\eta(\bE|X_{s\wedge\tau}-Y_{s\wedge\tau}|)\dif
s,
\de
which implies by Bihari inequality that for any $t>0$ and
$|x_0-y_0|<\eta$
$$
\bE|X_{t\wedge\tau}-Y_{t\wedge\tau}|\leq |x_0-y_0|^{\exp\{-\l_0t/2\}}
$$
and thus \be \bE(t\wedge\tau)\leq |x_0-y_0|^{1-\a} +\frac{\l_0
t}{2}\rho_\eta(|x_0-y_0|^{\exp\{-\l_0
t/2\}})\cdot|x_0-y_0|^{-\a}.\label{Op1} \ee

Taking
$\a=\exp\{-\l_0 T\}/2$, there exists an
$0<\eta'<\eta$ such that for any $|x_0-y_0|<\eta'$
\be\label{Eq9}
\bE((2T)\wedge\tau)\leq C_{T,\l_0,\eta'}\cdot|x_0-y_0|^{\exp\{-\l_0
T\}/2}.
\ee
But by Girsanov's theorem, $(\tilde W_t)_{t\in[0,T]}$ is still a
$n$-dimensional Brownian motion under the new probability measure
$U_T\cdot\bP$. Note that $(Y_t, \tilde{K}_t)$ also solves
$$
\dif Y_t+A(Y_t)\dif t\ni b(Y_t)\dif t+\sigma(Y_t)\dif \tilde W_t,\ \
\ Y_0=y_0.
$$
So, the law of $X_T(y_0)$ under $\bP$ is the same as that of
$Y_T(y_0)$ under $U_T\cdot\bP$.
Thus by (\ref{Eq7}), (\ref{Eq8}) and (\ref{Eq9}), for any $f\in
B_b(\mR^d)$,
\ce
&&|P_Tf(x_0)-P_Tf(y_0)|=|\bE(f(X_T(x_0))-U_T\cdot f(Y_T(y_0)))|\\
&\leq&\bE\left|(1-U_T)\cdot f(X_T(x_0))\cdot 1_{\{\tau\leq
T\}}\right|+\bE\left|(f(X_T(x_0))-U_T\cdot f(Y_T(y_0)))\cdot 1_{\{\tau>T\}}\right|\\
&\leq&\|f\|_0\cdot\bE|1-U_T|+\|f\|_0\cdot\bE\left[(1+U_T)1_{\{\tau>T\}}\right]\\
&\leq& C_{T,\l_0,\l_2,\eta}\cdot\|f\|_0\cdot|x_0-y_0|^{\exp\{-\l_0T\}/4}.\de

Step 2: Now we prove the proposition under {\bf (H3)}.
By the Markov property of the solution, we only need to prove that for every $f\in B_b(\mR^d)$,
$x\mapsto P_tf(x)$ is continuous on $D_r$ for all $t\leq t_p,$ $p>d$
where $p$ and $t_p$ are specified in Lemma \ref{lemma31}. Set
$$
c_0:=\|f\|_\infty
$$
and
$$
\tau:=\inf\left\{t>0:\sup_{x\in D_r}|X_t(x)|>N\right\}.
$$
Let $\e>0$ be given. For $t\leq t_p,$ by Lemma \ref{lemma31} and
Chebyshev inequality, there exists  $N>r$ such that \be \bP(\tau\leq
T)=\bP\left(\sup_{x\in D_r,t\leq t_p}|X_t(x)|>N\right)\leq \bE\left[\sup_{x\in D_r,
t\leq t_p}|X_t(x)|^p\right]/N^p<\e. \ee Define
$$
\tilde{\sigma}(x):=\sigma(x), ~~~\forall |x|\leq N.
$$
Extend $\tilde{\sigma}$ to the whole $\mR^d$ such that it satisfies
the condition {\bf (H1)} to {\bf (H3)}. Denote by $\tilde{X}_t(x)$ the solution
to (\ref{Eq2}) with $\sigma$ replaced by $\tilde{\sigma}$. By Step
1, there exists a $\delta>0$ such that if $|x-y|<\delta$ and $x,y\in
D_r$, \be |\bE[f(\tilde{X}_t(x))]-\bE[f(\tilde{X}_t(y))]|<\e. \ee
Hence \ce
&&|\bE[f(X_t(x))]-\bE[f(X_t(y))]|\\
&\leq &|\bE[(f(X_t(x))-f(X_t(y)))1_{(\tau>T)}]|+
|\bE[(f(X_t(x))-f(X_t(y)))1_{(\tau\leq T)}]|\\
&\leq &|\bE[(f(\tilde{X}_t(x))- f(\tilde{X}_t(y)))1_{(\tau>T)}]|+2c_0\e\\
&\leq &|\bE[f(\tilde{X}_t(x))- f(\tilde{X}_t(y))]|+
|\bE[(f(\tilde{X}_t(x))- f(\tilde{X}_t(y)))1_{(\tau\leq T)}]|+2c_0\e\\
&\leq & (1+4c_0)\e.\de
\end{proof}
Now we are in a position to complete the proof of Theorem \ref{Th3}.
\begin{proof}
$(i)$ By It\^o's formula and {\bf ($\mathbf{H}4$)}, we get \ce
\bE|X_t|^2&=&|x_0|^2+2\int_{0}^{t}\bE\<X_s,b(X_s)\>\dif
s-2\int_{0}^{t}\bE\<X_s,\dif
K_s\>+\int_{0}^{t}\bE\|\sigma(X_s)\|_\HS^2\dif
s\\&\leq&|x_0|^2+\int_{0}^{t}\bE(-\lambda_3|X_s|^p+\lambda_4)\dif
s.\de Taking derivatives with respect to $t$ and using H\"older's
inequality give \ce \frac{\dif \bE|X_t|^2}{\dif
t}&\leq&-\lambda_3\bE|X_t|^p+\lambda_4\leq-\lambda_3(\bE|X_t|^2)^{p/2}+\lambda_4.\de
Since $\lambda_3>0$ we have for all $t>0$,\ce
\frac{1}{t}\int_{0}^{t}\bE|X_s|^2\dif s\leq\l_4/\l_3.\de Therefore
by Krylov-Bogoliubov's method (see \cite{Ce}), there exists an
invariant probability measure $\mu$. As we have just proved, $P_t$
is strong Feller and irreducible.
Then, again by \cite{Ce}, $\mu$ is equivalent to each $P_t(x,\cdot)$ with $x\in\overline{D(A)},t>0$ and consequently $(i)$ holds. \\

$(ii)$ If $p>2,$ consider the following ODE:
$$
f'(x)=-\lambda_3f(x)^{p/2}+\lambda_4, ~~ f(0)=|x_0|^2.$$ By the
comparison theorem (cf. \cite{Ce}), there exists some $C>0$ such
that
$$
\bE|X_t|^2\leq f(t)\leq C(1+t^{2/(2-p)}). $$ We also have
$$
\inf_{x_0\in B(0,r)}P_t(x_0,B(0,a))>0,~~  \forall r,a>0,~t>0$$
because of the strong Feller property and irreducibility. Therefore
(ii) holds due to Theorem 2.5 (b) and Theorem 2.7 in \cite{Go-Ma}.
\end{proof}


\begin{thebibliography}{999}

\bibitem{ce1} C\'epa, E.:  \'Equations diff\'erentielles stochastiques
multivoques. Lect. Notes in Math. {\it S\'em. Prob. XXIX} (1995)
86-107.

\bibitem{ce3}C\'epa, E. et Jacquot, S.:  Ergodicit\'e D'in\'egalit\'es Variationnelles Stochastiques. {\it
Stochastics and Stochastics Reports.} Vol. 63  (1997) pp. 41-64.

\bibitem{Ce}Cerrai, S.: Second order PDE's in finite and infinite dimension. A probabilistic approach.
Lecture Notes in Mathematics, 1762. Springer-Verlag, Berlin, 2001.
x+330 pp.

\bibitem{Go-Ma}Goldys, B. and Maslowski, B.: Exponential ergodicity
for stochastic reaction-diffusion equations.  Stochastic partial
differential equations and applications---VII, 115--131, Lect. Notes
Pure Appl. Math., 245, Chapman  Hall/CRC, Boca Raton, FL, 2006.


\bibitem{rz1}Ren, J. and Zhang, X.: Stochastic flows for SDEs with non-Lipschitz coefficients,
{\it Bull. Sci. Math.}, Vol. 127 (2003), 739-754.

\bibitem{Zh}Zhang, X.: Exponential ergodicity of non-Lipschitz
stochastic differential equations. {\it Proc. Amer. Math. Soc.} 137  (2009), 329-337.

\end{thebibliography}
\end{document}